\documentclass[12pt]{amsart}
\usepackage{amssymb}
\usepackage{textcomp}
\usepackage[all]{xy}

\title{Almost non-degenerate abelian fibrations}
\author{Ying Zong}

\address{Department of Mathematics\\University of Toronto}
\email{zongying@math.utoronto.ca}
\date{} 

\begin{document}

\maketitle


1. \emph{Introduction.}
\smallskip

Let $S$ be an algebraic space and $X$ an $S$-algebraic space. One says that $X$ is a \emph{non-degenerate} $S$-\emph{abelian fibration} if there is an $S$-abelian algebraic space $A$ such that $X$ is an $A$-torsor on $S$ for the \'{e}tale topology in which case $A$ is then the \emph{albanese} of this fibration and hence is uniquely determined. In the following \S 3 we define the notion for $X$ to be an \emph{almost non-degenerate} $S$-\emph{abelian fibration} and with each such associate its \emph{albanese} over its \emph{ramification} $S$-\emph{stack} and our goal is to support the \emph{Principle} : \emph{Non-uniruled abelian fibrations are almost non-degenerate}.
This support is however very partial, as the purity theorem below (\cite{grothendieck_abelian} 4.5) on which our arguments at one step crucially rely fails in positive or mixed characteristics :
\smallskip

\emph{If $T$ is a locally noetherian regular algebraic space of residue characteristics zero and if $U$ is an open sub-algebraic space of $T$ such that $\mathrm{codim}(T-U, T)\geq 2$, then the functor $A\mapsto A|U$, from the category of $T$-abelian algebraic spaces to the category of $U$-abelian algebraic spaces, is an equivalence.}

\smallskip

When the conclusion of this purity theorem holds, it suffices largely to consider non-uniruled abelian degenerations over spectra of discrete valuation rings, for which our results are in \S 4. The extension to the global situation in characteristic zero with hypothesis of purity and para-factoriality is in the last section \S 5. Besides some terminologies, in \S 2 there is a write-up of Gabber's theorem of \emph{purity of branch locus}. 

\smallskip

2. \emph{Rational curves, uniruled irreducible components, regular minimal models and purity of branch locus.}
\smallskip

Recall some terminologies.
\smallskip

Let $k$ be an algebraically closed field.
\smallskip

A $k$-\emph{rational curve} is by definition a separated integral $k$-scheme  $C$ of finite type such that $k(\eta)$ is a purely transcendental extension of $k$ of transcendence degree $1$, where $\eta$ denotes the generic point of $C$.
\smallskip

Let $V$ be a $k$-algebraic space. One says that $V$ \emph{does not contain $k$-rational curves} if every $k$-morphism from a $k$-rational curve $C$ to $V$ factors as $C\to \mathrm{Spec}(k)\to V$ for a certain $k$-point $\mathrm{Spec}(k)\to V$ of $V$.
\smallskip

If an algebraic space $V$ is quasi-separated locally of finite type over $k$ and does not contain $k$-rational curves, then the base change of $V$ to every algebraically closed extension $k'$ of $k$ does not contain $k'$-rational curves.
\smallskip

For a proper $k$-algebraic space $V$, saying that $V$ does not contain $k$-rational curves amounts to saying that every $k$-morphism from the $k$-projective line to $V$ factors through a $k$-point of $V$. Thus, as the $k$-projective line is simply connected, if a proper $k$-algebraic space does not contain $k$-rational curves, neither does its quotient by any finite \'{e}tale $k$-equivalence relation.
\smallskip

The following result of Murre and Chow, which has origin in one of Zariski's proofs of his Connectedness/Main Theorem, explains the significance of non-existence of rational curves.
\smallskip

{\bf Lemma 2.1.} --- \emph{Let $S$ be an algebraic space and $T\to S$ a morphism to $S$ from a connected locally noetherian regular algebraic space $T$. Let $X$ be a proper $S$-algebraic space whose geometric $S$-fibers do not contain rational curves.}
\smallskip

\emph{Then every $S$-morphism from a non-empty open sub-algebraic space of $T$ to $X$ extends uniquely to an $S$-morphism from $T$ to $X$.}

\begin{proof} This is \cite{grothendieck_rational} 2.6.

\end{proof}

\smallskip

Let $F$ be a finitely generated field over $k$. One says that $F$ is \emph{ruled} over $k$ or that $F/k$ is \emph{ruled}, if $F$ has a $k$-sub-field $K$ such that $\mathrm{Spec}(F)\to\mathrm{Spec}(K)$ is geometrically regular, that is, the extension $F/K$ is \emph{regular} in the sense of Weil, and furthermore such that $\overline{K}\otimes_KF$ is a purely transcendental extension of $\overline{K}$, where $\overline{K}$ is an algebraic closure of $K$. One says that $F$ is \emph{uniruled} over $k$ or that $F/k$ is \emph{uniruled}, if $F$ has a finite extension $F'/k$ which is ruled.
\smallskip

Let $V$ be a quasi-separated algebraic space locally of finite type over $k$, $\eta$ a maximal point of $V$ and $Z$ the irreducible component of $V$ with reduced induced structure with generic point $\eta$. If $k(\eta)/k$ is ruled (resp. uniruled), one says that $V$ is \emph{ruled} (resp. \emph{uniruled}) at $\eta$ and that $Z$ is a \emph{ruled} (resp. \emph{uniruled}) irreducible component of $V$.

\smallskip

{\bf Lemma 2.2.} --- \emph{An abelian variety over an algebraically closed field $k$ does not contain $k$-rational curves and thus in particular is not uniruled. A connected smooth algebraic group over an algebraically closed field is not ruled only if it is an abelian variety.}

\begin{proof} See \cite{neron_model} 9.2/4, 9.2/1.

\end{proof}

\smallskip

Let $S$ be the spectra of a discrete valuation ring and $t$ the generic point of $S$. Let $A$ be an $S$-abelian scheme and $X_t$ an $A_t$-torsor on $t$ for the \'{e}tale topology. Recall (\cite{raynaud_minimal}, p. 82, line -2) that $X_t$ with its $A_t$-action admits a unique extension to an $S$-scheme $X$ with an action by $A$ such that $X$ is projective flat over $S$, regular and such that the morphism
\[A\times_SX\to X\times_SX,\ (a, x)\mapsto (a+x, x)\] is finite surjective. The geometric $S$-fibers of $X$ are irreducible and do not contain rational curves. Following N\'{e}ron--Raynaud, one says that $X$ is the \emph{regular} $S$-\emph{minimal model of} $X_t$. The formation of regular minimal models commutes with every formally smooth faithfully flat base change $S'\to S$ of spectra of discrete valuation rings. 
\smallskip

{\bf Lemma 2.3.} --- \emph{Let $S$ be the spectra of a discrete valuation ring, $t$ the generic point of $S$, $A$ an $S$-abelian scheme and $X_t$ an $A_t$-torsor on $t$ for the \'{e}tale topology. Assume that $X_t$ extends to an $S$-algebraic space $X$ which is proper over $S$, regular and with an irreducible closed $S$-fiber.}
\smallskip

\emph{Then $X$ is the regular $S$-minimal model of $X_t$.}

\begin{proof} Let $s$ be the closed point of $S$ and $x$ the generic point of $X_s$. Notice that $X$ is connected and that there is an open neighborhood of $x$ in $X$ which is a scheme (\cite{raynaud_specialization} 3.3.2).
Let $Y$ be the regular $S$-minimal model of $X_t$. The identity $X_t=Y_t$ extends uniquely by (2.1) to an $S$-morphism $f: X\to Y$. Being proper birational, $f$ maps $x$ to the generic point of $Y_s$ and induces an isomorphism between spectra of discrete valuation rings
\[\mathrm{Spec}(\mathcal{O}_{X, x})\ \widetilde{\to}\ \mathrm{Spec}(\mathcal{O}_{Y, f(x)}).\] By the theorem of \emph{purity of branch locus} (2.4) below, $f$ is \'{e}tale and hence is an isomorphism.

\end{proof}

\smallskip

The following theorem answers the question EGA IV 21.12.14 (v).
\smallskip

{\bf Theorem 2.4}  (Gabber). --- \emph{Let $f: X\to Y$ be a morphism, essentially of finite type, from a normal scheme $X$ to a locally noetherian regular scheme $Y$ such that $f$ is essentially \'{e}tale at all points of $X$ of codimension $\leq 1$.}
\smallskip

\emph{Then $f$ is essentially \'{e}tale.}

\begin{proof} When $f$ is moreover finite, $f$ is \'{e}tale by SGA 2 X 3.4. It follows that $f$ is essentially \'{e}tale at a point $x\in X$ if it is essentially quasi-finite at $x$. For, letting $X_{(x)}$ (resp. $Y_{(f(x))}$, resp. $f_{(x)}$) be the henselization of $X$ (resp. $Y$, resp. $f$) at $x$ (resp. $f(x)$, resp. $x$), one can apply \emph{loc.cit.} to $f_{(x)}: X_{(x)}\to Y_{(f(x))}$, as $f_{(x)}$ is finite by Zariski's Main Theorem. 
\smallskip

So it amounts to showing that $f$ is essentially quasi-finite. It suffices to see that $f$ is essentially quasi-finite at a point $x\in X$ if it is at all generalizations of $x$.
\smallskip
 
One may assume that $X$ is local of dimension $\geq 2$ with closed point $x$, that $Y$ is local with closed point $f(x)$, that $f|(X-\{x\})$ is essentially \'{e}tale and that $x=f^{-1}(f(x))$. Notice that $f$ is essentially quasi-finite at $x$ if and only if the extension $k(x)/k(f(x))$ is finite.
\smallskip

--- \emph{Reduction to the case where $Y$ is excellent }:
\smallskip

The completion $Y'$ of $Y$ along $f(x)$ is excellent, so is $X'=X\times_YY'$, and the normalization $X''$ of the excellent scheme $X'$ is finite over $X'$. Write $y'$ for the closed point of $Y'$ and let $x'$ be the unique point of $X'$ with image $x$ in $X$ and with image $y'$ in $Y'$. The projection $X''\to X'$ restricts to an isomorphism over $X'-\{x'\}$, since $X'-\{x'\}$, being essentially \'{e}tale over $Y'$, is regular. 
It suffices to show that the composition $X''\to X'\to Y'$ is essentially quasi-finite at a point $x''$ of $X''$ above $x'$. For then the extension $k(x'')/k(y')$, a priori $k(x')/k(y')$ as well, is finite and hence the extension $k(x)/k(f(x))$ is finite.
\smallskip

--- \emph{Assume $Y$ excellent. Reduction to the case $\mathrm{dim}(X)=2$ }:
\smallskip

One has $\mathrm{dim}(Y)>0$. Let $h\in\Gamma(Y, \mathcal{O}_Y)$ be part of a regular system of parameters at $f(x)$, let $Y'=V(h)$ and let $X'=V(h\mathcal{O}_X)=X\times_YY'$. The normalization $X''$ of the excellent scheme $X'$ is finite over $X'$ and the projection $X''\to X'$ restricts to an isomorphism over $X'-\{x\}$, as $X'-\{x\}$, being essentially \'{e}tale over $Y'$, is regular. It suffices to show that the composition $X''\to X'\to Y'$ is essentially quasi-finite at a point $x''\in X''$ above $x\in X'$, for the extensions $k(x'')/k(f(x))$ and $k(x)/k(f(x))$ are then finite. Note that $X''$ is purely of dimension $\mathrm{dim}(X)-1$, which is $\geq 2$ if and only if $\mathrm{dim}(X)>2$.
\smallskip

--- \emph{Assume $Y$ excellent and $\mathrm{dim}(X)=2$. Then $f^{!}\mathcal{O}_Y=\mathcal{O}_X$ }:
\smallskip

Now $Y$ being regular local, $\mathcal{O}_Y$ is a dualizing object on $Y$, so is $f^{!}\mathcal{O}_Y$ on $X$. As $X$ being normal local of dimension $2$ is Cohen--Macaulay, $f^{!}\mathcal{O}_Y$ has one only non-zero cohomology, say $L$, which is Cohen--Macaulay, concentrated in degree, say $d$. Let $U=X-\{x\}$ and $j: U\to X$ the open immersion. Since $f$ is essentially \'{e}tale on $U$, 
\[(fj)^{!}\mathcal{O}_Y=j^{!}f^{!}\mathcal{O}_Y=j^*f^{!}\mathcal{O}_Y=j^*(L[-d])=(L|U)[-d]\] is canonically isomorphic to $\mathcal{O}_U$ in $D(U, \mathcal{O}_U)$. That is, $d=0$ and $L|U$ is canonically isomorphic to $\mathcal{O}_U$. Such an isomorphism uniquely extends to an isomorphism from $\mathcal{O}_X$ onto $L=f^{!}\mathcal{O}_Y$, as $\mathrm{prof}_x(\mathcal{O}_X)=\mathrm{prof}_x(L)=2$.
\smallskip

--- \emph{Assume $Y$ excellent and $\mathrm{dim}(X)=2$. Then $X$ is regular }:
\smallskip

One has $\mathrm{dim}(Y)>0$. Let $h\in\Gamma(Y, \mathcal{O}_Y)$ be part of a regular system of parameters at $f(x)$, let $Y'=V(h)$, $i: Y'\to Y$ the canonical closed immersion, $f'=f\times_YY': X'\to Y'$, let the normalization of $X'=V(h\mathcal{O}_X)=X\times_YY'$ be $p: X''\to X'$ and let $f''=f'p: X''\to X'\to Y'$ be the composition. 
\smallskip

Notice that $f: X\to Y$ and $i: Y'\to Y$ are tor-independent over $Y$. So, from the identity $f^{!}\mathcal{O}_Y=\mathcal{O}_X$, it follows that $f^{'!}\mathcal{O}_{Y'}=\mathcal{O}_{X'}$. 
\smallskip

Notice next that $f''$ is essentially of complete intersection, for $Y'$ is regular and $X''$, being normal of dimension $1$, is regular. As $f''$ is furthermore essentially \'{e}tale at all points of $X''$ above $X'-\{x\}$, the object $f^{''!}\mathcal{O}_{Y'}$ has a unique non-zero cohomology in degree $0$, which is an invertible $\mathcal{O}_{X''}$-module, and there is a canonical homomorphism, the fundamental class of $f''$, $c(f''): \mathcal{O}_{X''}\to f^{''!}\mathcal{O}_{Y'}$. 
This $c(f'')$ corresponds by duality to a morphism
$p_{*}\mathcal{O}_{X''}\to f^{'!}\mathcal{O}_{Y'}$, which is injective and which when composed with the canonical homomorphism
$\mathcal{O}_{X'}\to p_{*}\mathcal{O}_{X''}$ induces the above identity $\mathcal{O}_{X'}=f^{'!}\mathcal{O}_{Y'}$, as one verifies at each point of $X'-\{x\}$ over which $p$ is an isomorphism. 
\smallskip

So $\mathcal{O}_{X'}=p_{*}\mathcal{O}_{X''}$ and $X'$ is regular. So $X$ is regular.
\smallskip

--- \emph{Assume $Y$ excellent, $X$ regular and $\mathrm{dim}(X)=2$. Then $f$ is essentially \'{e}tale }:
\smallskip

Now $f$ is essentially of complete intersection. Its cotangent complex $L_f$ is a perfect complex in $D^{[-1,0]}(X, \mathcal{O}_X)$. With $L_f$ one associates a canonical ``theta divisor'' homomorphism
\[\theta: \mathcal{O}_X\to \mathrm{det}(L_f),\] which, as $f|(X-\{x\})$ is essentially \'{e}tale, is an isomorphism outside $x$. As $\mathrm{prof}_x(X)=2$ and $\mathrm{det}(L_f)$ is an invertible $\mathcal{O}_X$-module, $\theta$ is an isomorphism also at $x$. So $f$ is essentially \'{e}tale at $x$ by the Jacobian Criterion.

\end{proof}

\smallskip

3. \emph{Almost non-degenerate abelian fibrations.}
\smallskip

{\bf Definition 3.1.} --- \emph{Let $S$ be an algebraic space. We say that an $S$-algebraic space $X$ is an almost non-degenerate $S$-abelian fibration if there exists in the category of $S$-algebraic spaces a groupoid whose nerve $(X_., d_., s_.)$ satisfies the conditions a) and b) }:
\smallskip

\emph{a) $X=X_o$.}
\smallskip

\emph{b) $d_1: X_1\to X_o$ is the structural morphism of an $X_o$-abelian algebraic space with zero section $s_o: X_o\to X_1$.}
\smallskip

\emph{If for every smooth morphism $S'\to S$ $\mathrm{Coker}(d_o\times_SS', d_1\times_SS')=S'$ in the category of $S'$-algebraic spaces, we say that $X$ is a geometric almost non-degenerate $S$-abelian fibration. If $\mathrm{Coker}(d_o\times_SS', d_1\times_SS')=S'$ holds in the category of $S'$-algebraic spaces for every base change $S'\to S$, we say that $X$ is a universal almost non-degenerate $S$-abelian fibration.}
\smallskip

\emph{We say that the $S$-stack $[X_.]$ defined by the groupoid $X_.$ is the ramification $S$-stack of $X$.}
\smallskip

\emph{By a morphism $X\to X'$ of almost non-degenerate $S$-abelian fibrations, we mean a morphism of $S$-groupoids $X_.\to X'_.$.}
\smallskip

Given every such morphism $f$ from $X$ to $X'$, $f_1: X_1\to X'_1$ is an $f_o$-homomorphism with respect to the abelian algebraic space structures $(d_1, s_o, +)$ on $X_1/X_o$ and on $X'_1/X'_o$ (``Geometric Invariant Theory'' 6.4).
\smallskip

The isomorphism (``renversement des fl\`{e}ches'') $X_1\to X_1$ transports to $d_o: X_1\to X_o$ from $d_1: X_1\to X_o$ an $X_o$-abelian algebraic space structure with zero section $s_o: X_o\to X_1$.
\smallskip

The structural morphism $X\to S$ factors, which is the essence of (3.1), as
\[X\to [X_.]\to S\] in such a tautological way that the projection $X\to [X_.]$ is a torsor for the lisse-\'{e}tale topology under an $[X_.]$-abelian algebraic stack $A$ verifying $A\times_{[X_.]}X=X_1$. We call $A/[X_.]$ the \emph{albanese} of $X/S$, cf. \cite{basic} 7.2+7.3.
\smallskip

Notice that $X/S$ is flat (resp. proper, resp. of finite presentation) if and only if $[X_.]/S$ is flat (resp. proper, resp. of finite presentation).
\smallskip

For every $S$-algebraic space $S'$, $X\times_SS'$ is an almost non-degenerate $S'$-abelian fibration with defining $S'$-groupoid $X_.\times_SS'$.
\smallskip

Let $E$ be an algebraic stack. We say that an $E$-algebraic stack $F$ is an \emph{almost non-degenerate} $E$-\emph{abelian fibration} if there is a smooth surjective morphism from an algebraic space $S$ to $E$ such that $X=F\times_ES$ is an almost non-degenerate $S$-abelian fibration endowed with a descent datum on its $S$-groupoid $X_.$ relative to $S\to E$. We say that $F/E$ is \emph{geometric} (resp. \emph{universal}) if $X/S$ is.
\smallskip

If the $S$-groupoid $X_.$ is simply connected with $\mathrm{Coker}(d_o, d_1)=S$, that is, if $[X_.]=S$, then $X/S$ is a \emph{non-degenerate abelian fibration}, namely, a torsor for the \'{e}tale topology under an $S$-abelian algebraic space. One says correspondingly that an almost non-degenerate abelian fibration $F/E$ as above is \emph{non-degenerate} when $F\times_ES/S$ is non-degenerate. 
\smallskip

4. \emph{Non-uniruled abelian fibrations over spectra of discrete valuation rings.}
\smallskip

{\bf Theorem 4.1.} --- \emph{Let $S$ be the spectra of a discrete valuation ring, $t$ the generic point of $S$, $s$ the closed point and $\overline{s}$ the spectrum of an algebraic closure of $k(s)$. Let $A_t$ be a $t$-abelian variety and $X_t$ an $A_t$-torsor on $t$ for the \'{e}tale topology. Assume that $X_t$ extends to a separated faithfully flat $S$-algebraic space $X$ of finite type such that not all irreducible components of $X_{\overline{s}}$ are uniruled.}
\smallskip

\emph{Then }:
\smallskip

\emph{a) There is a spectrum $t'$ of a finite separable extension of $k(t)$ such that, if $S'$ denotes the normalization of $S$ in $t'$, $A_t\times_tt'$ extends to an $S'$-abelian scheme $A'$.}
\smallskip

\emph{b) Exactly one irreducible component of $X_{\overline{s}}$ is not uniruled. In particular, $X_{\overline{s}}$ is irreducible if it does not have uniruled irreducible components.}
\smallskip

{\bf Theorem 4.2.} --- \emph{Keep the assumptions of $(4.1)$. Assume furthermore that $X_{\overline{s}}$ is connected, proper and separable.}

\smallskip

\emph{Then $X$ is proper over $S$ and normal, $X_{\overline{s}}$ has a unique non-uniruled irreducible component, $A_t$ extends to an $S$-abelian scheme $A$, the regular $S$-minimal model of $X_t$ is an $A$-torsor $F$ on $S$ for the \'{e}tale topology, the canonical birational map from $X$ to $F$ is an $S$-rational map $p$ whose domain of definition contains all points where $X$ is geometrically factorial of equal characteristic or regular, $p$ is \'{e}tale at precisely the points of $\mathrm{Dom}(p)$ outside the image in $X$ of the uniruled irreducible components of $X_{\overline{s}}$ and $p^{-1}$ extends to a proper $S$-birational morphism from $F$ onto $X$ if $X_{\overline{s}}$ does not contain $\overline{s}$-rational curves.}

\smallskip

{\bf Corollary 4.3.} --- \emph{Keep the assumptions of $(4.1)$. Assume furthermore that $X_{\overline{s}}$ is proper, separable and does not have uniruled irreducible components and that $X$ is at each of its geometric codimension $\geq 2$ points either geometrically factorial of equal characteristic or regular.}
\smallskip

\emph{Then $A_t$ extends to an $S$-abelian scheme $A$ and $X$ is an $A$-torsor on $S$ for the \'{e}tale topology.}

\smallskip

{\bf Lemma 4.4.} --- \emph{Let $S$ be the spectra of an excellent discrete valuation ring with generic point $t$. Let $S'\to S$ be a surjective morphism to $S$, essentially of finite type, from a local integral scheme $S'$ of dimension $1$ with generic point $t'$ such that the extension $k(t')/k(t)$ is regular of transcendence degree $d\geq 1$.}
\smallskip

\emph{Then there is an $S$-scheme $S_o$, which is the spectra of a discrete valuation ring and which is quasi-finite surjective over $S$, and there exists a sequence of morphisms, $S_d\to S_{d-1}\to\cdots\to S_1\to S_o$, each of which is smooth, surjective, purely of relative dimension $1$, with geometrically connected fibers, such that the normalization $S''$ of $S'\times_SS_o$ is $S_o$-isomorphic to a localization of $S_d$. In particular, $S''$ is essentially smooth over $S_o$.}

\begin{proof} This is \cite{dejong} 2.13.

\end{proof}

\smallskip

{\bf Lemma 4.5.} --- \emph{Let $p: V\to S$ be a morphism, essentially of finite type, from a regular local scheme $V$ of dimension $1$ with closed point $v$ to a regular local scheme $S$ of dimension $>1$ with closed point $s$ such that $p$ is birational and that $p(v)=s$.}
\smallskip

\emph{Then $k(v)$ has a sub-$k(s)$-extension $K$ such that the extension $k(v)/K$ is purely transcendental.}

\begin{proof} Recall the argument of Zariski. Let $f: X\to S$ be the blow up of $S$ along $s$ (EGA II 8.1.3). Then $X$ is regular (EGA IV 19.4.3, 19.4.4), $f$ is proper and there is one only $S$-morphism $p_1: V\to X$ by the valuative criterion of properness. Write $s_1=p_1(v)$, $S_1=\mathrm{Spec}(\mathcal{O}_{X, s_1})$. Denote the canonical morphism $V\to S_1$ again by $p_1$. Blowing up $S_1$ along $s_1$ and localizing, one obtains similarly $p_2: V\to S_2$. Continuing this way, one finds a projective system of regular local schemes ``$\underleftarrow{\mathrm{Lim}}$'' $S_i$, indexed by $\mathbf{N}$, with $S_o=S$, with each transition morphism being birational, local and essentially of finite type. There is a unique $S$-morphism $(p_i): V\to$ ``$\underleftarrow{\mathrm{Lim}}$'' $S_i$ with $p_o=p$. Write $s_n$ for the closed point of $S_n$, $n\in\mathbf{N}$.
\smallskip

It suffices to show that the projective system ``$\underleftarrow{\mathrm{Lim}}$'' $S_i$ is essentially constant and has $V$ as its limit. For, if $n$ is the smallest integer such that $V=S_n$, the extension $k(v)/k(s_{n-1})$ is purely transcendental.
\smallskip

Let $\mathfrak{m}$ denote the ideal of $\mathcal{O}_S$ defining the closed point $s$. For every coherent $S$-ideal $I$ that is non-zero and distinct from $\mathcal{O}_S$, $I\mathcal{O}_X$ is a non-zero sub-ideal of $\mathfrak{m}\mathcal{O}_X=\mathcal{O}_X(1)$ (EGA II 8.1.7). Thus 
\[I\mathcal{O}_X\otimes_X\mathcal{O}_X(-1)=I_X(-1)\] is a non-zero ideal of $\mathcal{O}_X$. Write $I_1$ for the localization of $I_X(-1)$ at $s_1$. Then $I\mathcal{O}_V$ is \emph{strictly} contained in $I_1\mathcal{O}_V$. If $I_1$ is distinct from $\mathcal{O}_{S_1}$ and if $S_2$ is distinct from $S_1$, one obtains as above an ideal $I_2$ of $\mathcal{O}_{S_2}$ such that $I_1\mathcal{O}_V$ is strictly contained in $I_2\mathcal{O}_V$. As $V$ is noetherian, this process of producing ideals $I_1, I_2, \cdots$ from a given $I$ eventually stops. 
That is, for a certain integer $N\geq 1$, either $S_N$ is of dimension $1$ (thus is equal to $V$) or $I_N=\mathcal{O}_{S_N}$. Assume that the latter case holds for each coherent $S$-ideal $I$ that is distinct from $0$ and $\mathcal{O}_S$. For otherwise the above assertion is already proven. 
\smallskip

Choose $a_j, b_j\in\Gamma(S, \mathcal{O}_S)$, $a_j\mathcal{O}_V=b_j\mathcal{O}_V$ non-zero, $j=1,\cdots, d$, where $d=\mathrm{deg.tr.}(k(v)/k(s))$, such that the images of the fractions
\[a_1/b_1,\ \cdots,\ a_d/b_d\] in $k(v)$ form a basis of transcendence over $k(s)$. 
By the assumption just made applied to the ideals $a_j\mathcal{O}_S, b_j\mathcal{O}_S$, there exists an integer $N_j\geq 1$ for each $j=1,\cdots, d$ such that either $a_j/b_j$ or $b_j/a_j$ is a section of $\mathcal{O}_{S_{N_j}}$ over $S_{N_j}$. Being invertible on $V$, $a_j/b_j$ or equivalently $b_j/a_j$ is invertible on $S_{N_j}$. 
Hence, for an integer $N\geq 1$, $k(v)$ is algebraic over $k(s_N)$ and $V$ is essentially quasi-finite over $S_N$. By Zariski's Main Theorem, $V$ is isomorphic to $S_N$.

\end{proof}

\smallskip

4.6. \emph{Proof of} (4.1).
\smallskip

See \cite{basic} \S 10 for an application of (4.1) \emph{b}).
\smallskip

We reduce the proof to that of (4.2).
\smallskip

--- \emph{Reduction to the case where $S$ is strictly henselian }:
\smallskip

If $\overline{t}$ denotes the spectrum of a separable closure of $k(t)$ and if $\ell$ is a prime number prime to the characteristic of $k(s)$, the claim \emph{a}) is equivalent to the claim that the $\ell$-adic monodromy representation associated with the $t$-abelian variety $A_t$,
\[\rho_{\ell, \overline{t}}: \pi_1(t, \overline{t})\to \mathrm{GL}(H^1(A_{\overline{t}}, \mathbf{Q}_{\ell})),\] when restricted to an inertia subgroup relative to $S$, has finite image (\cite{neron_model} 7.4/5). For both assertions \emph{a}), \emph{b}), replacing $S$ by its strict henselization $S_{(\overline{s})}$ at $\overline{s}$, $X$ by $X\times_SS_{(\overline{s})}$ and $A_t$ by $A_t\times_tt^{hs}$, where $t^{hs}$ denotes the generic point of $S_{(\overline{s})}$, one may assume that $S$ is strictly henselian.
\smallskip

--- \emph{Assume $S$ strictly henselian. Reduction to the case where $S$ is complete }:
\smallskip

Let $S'$ be the completion of $S$ along $s$, $t'$ the generic point of $S'$ and $\overline{t'}$ the spectrum of a separable closure of $k(t')$ containing a separable closure $k(\overline{t})$ of $k(t)$. 
\smallskip

The projection $t'\to t$ induces by SGA 4 X 2.2.1 an isomorphism 
\[\pi_1(t', \overline{t'})\ \widetilde{\to}\ \pi_1(t, \overline{t}),\] relative to which the proper base change isomorphism (SGA 4 XII 5.1)
\[H^1(A_{\overline{t}}, \mathbf{Q}_{\ell})\ \widetilde{\to}\ H^1(A_{\overline{t'}}, \mathbf{Q}_{\ell})\] is equivariant. Replacing $S$ by $S'$, $X$ by $X\times_SS'$ and $A_t$ by $A_t\times_tt'$, one may assume moreover that $S$ is complete.
\smallskip

--- \emph{Assume $S$ strictly henselian and complete. Reduction to the case where $X$ is a proper flat $S$-scheme }:
\smallskip

Let the maximal points of $X_s$ be $x_1,\cdots, x_n$. There is an open neighborhood $V$ of $\{x_1,\cdots, x_n\}$ in $X$ which is a scheme (\cite{raynaud_specialization} 3.3.2). Replacing $X$ by $V$, one may assume that $X$ is a separated flat $S$-scheme of finite type. By Nagata, $X$ admits an open $S$-immersion into a proper $S$-scheme $P$. One may then replace $X$ by its closed image in $P$ and assume that $X$ is a proper flat $S$-scheme. 
\smallskip

--- \emph{Assume $S$ strictly henselian, complete and that $X$ is a proper flat $S$-scheme. Reduction to the case where $X_s$ is separable over $k(s)$ }: 
\smallskip

Notice that $X$ is integral. For, one has $\mathrm{Ass}(X)=\mathrm{Ass}(X_t)$, as $X$ is flat over $S$ (EGA IV 3.3.1). Notice also that $S$, being complete, is excellent. Let the maximal points of $X_s$ be $x_1,\cdots, x_n$. By applying (4.4) to the morphisms $\mathrm{Spec}(\mathcal{O}_{X, x_i})\to S$, one finds an $S$-scheme $S_o$, which is the spectra of a discrete valuation ring and which is finite surjective over $S$, such that the normalization of $X\times_SS_o$, $X_o$, is smooth over $S_o$ at all  maximal points of the closed fiber of $X_o/S_o$. Again, at least one irreducible component of the geometric closed fiber of $X_o/S_o$ is not uniruled. Observe that $X/S$ verifies the claims \emph{a})+\emph{b}) as long as $X_o/S_o$ does. Thus, replacing $S$ by $S_o$, $X$ by $X_o$ and $A_t$ by $A_t\times_tt_o$, where $t_o$ denotes the generic point of $S_o$, one may assume that $X$ is smooth over $S$ at all maximal points of $X_s$.
\smallskip

The fiber $X_{\overline{s}}$ is connected, as $X$ is proper over $S$. To finish it suffices to apply (4.2).
\smallskip

\begin{flushright}
$\square$
\end{flushright}

\smallskip

4.7. \emph{Proof of} (4.2), (4.3).
\smallskip

See \cite{basic} 10.2 for an application of (4.3).
\smallskip

One may assume that $S$ is strictly henselian.
\smallskip

Being separated faithfully flat over $S$ of finite type with geometric fibers proper and connected, $X$ is by EGA IV 15.7.10 proper over $S$. Moreover,  as $X_s$ is separable, $X$ is normal. And, the open sub-algebraic space $V$ of $X$ consisting of all points at which $X\to S$ is smooth is $S$-schematically dense in $X$. 
\smallskip

Fix an $S$-section $o$ of $V$ (EGA IV 17.16.3) by means of which one identifies $X_t=V_t$ with $A_t$. 
\smallskip

Consider a smoothening of $X$ (\cite{neron_model} 3.1/1, 3.1/3), $f: X'\to X$, which is so constructed as the composition of a finite sequence of blow-ups with centers lying above the complement of $V$ in $X$. Let $W'$ be the open sub-algebraic space of $X'$ consisting of all points at which $X'\to S$ is smooth. Naturally, $V$ can be identified with an open sub-algebraic space of $W'$ and one has $W'_t=V_t$. Let $d=\mathrm{dim}(A_t)$. Choose a non-zero section $\omega'\in\Gamma(W', \Omega^d_{W'/S})$ such that the support of the divisor $\mathrm{Div}_{W'}(\omega')$ is \emph{strictly} contained in $W'_s$. Such a section $\omega'$ clearly exists and is unique up to multiplication by a unit of $\Gamma(S, \mathcal{O}_S)$. 
\smallskip

The morphism $W'_t=A_t$, where the identification is provided by the above section $o$, has a unique extension to an $S$-morphism
\[p: W'\to A,\] where $A$ denotes the $S$-N\'{e}ron model of $A_t$. 
\smallskip

In the language of N\'{e}ron models, $W'$ is a weak $S$-N\'{e}ron model of $W'_t$ (\cite{neron_model} 3.5/1, 3.5/2). Its open sub-algebraic space $U'=W'-\mathrm{Supp}(\mathrm{Div}_{W'}(\omega'))$ admits a canonical $S$-birational group law (\emph{loc.cit.} 4.3/5) which extends the group structure of $U'_t=A_t$ over $t$. And, the restriction of $p$ to $U'$, $p|U': U'\to A$, is an $S$-schematically dense open immersion, which solves the universal problem of extending the $S$-birational group law of $U'$ to an $S$-group law (\emph{loc.cit.} 4.3/6, 5.1/5). 
\smallskip

Let the maximal points of $X_s$ be $x_1,\cdots, x_n$. They lie in $V$ and can thus be considered as points of $W'$. Notice that there is an open neighborhood of $\{x_1,\cdots, x_n\}$ in $W'$ which is a scheme (\cite{raynaud_specialization} 3.3.2).
\smallskip

Observe that if a point $x$ among $x_1,\cdots, x_n$ belongs to $V-U'$, that is, if $\mathrm{Div}_{W'}(\omega')_{x}$ is not zero, then 
\[p: \mathrm{Spec}(\mathcal{O}_{W', x})\to \mathrm{Spec}(\mathcal{O}_{A, p(x)})\] is not an isomorphism. This implies by (4.5) that $X_{\overline{s}}$ is uniruled at $x$. Here we have denoted again by $x$ the unique point of $X_{\overline{s}}$ that projects to $x$ in $X_s$. 
\smallskip

By hypothesis there is at least one point of $\{x_1,\cdots, x_n\}$, say $x_1$, at which $X_{\overline{s}}$ is not uniruled. One finds thus $x_1\in V\cap U'$ and that $A_{\overline{s}}$ is not uniruled at $p(x_1)$. So $A_{\overline{s}}$ is an $\overline{s}$-abelian variety (2.2). So $A$ is an $S$-abelian scheme, $U'_{\overline{s}}$ is irreducible, $X_{\overline{s}}$ is by (4.5) uniruled at all its maximal points other than $x_1$ and $p$ is \'{e}tale at precisely the points of $U'$. 
\smallskip

If $X_{\overline{s}}$ does not contain $\overline{s}$-rational curves, the rational map $p^{-1}$ extends by (2.1) to an $S$-morphism, hence a proper $S$-birational morphism, from $A$ onto $X$. 
\smallskip

Finally, being a trivial $A_t$-torsor, $X_t$ has a trivial $A$-torsor $F$ as its regular $S$-minimal model. The assertion on the points where the canonical birational map from $X$ to $F$ is \'{e}tale (resp. is  defined) follows by (2.4) (resp. (5.2)+(5.3) below). Recall (EGA IV 21.13.9, 21.13.11) that a noetherian normal local scheme of dimension $\geq 2$ is factorial (resp. geometrically factorial) if and only if it is para-factorial (resp. geometrically para-factorial) at all of its points (resp. geometric points) of codimension $\geq 2$. With the conditions of (4.3), such a birational map is everywhere defined and is an $S$-isomorphism.

\begin{flushright}
$\square$
\end{flushright}

\smallskip

{\bf Lemma 4.8.} --- \emph{Let $n, \delta_1,\cdots, \delta_n$ be integers $\geq 1$ such that the greatest common divisor of $\delta_1,\cdots, \delta_n$ is $1$. Let $S$ be an algebraic space, $A$ an $S$-abelian algebraic space, $X$ an $A$-torsor on $S$ for the \'{e}tale topology and $S_i$ an $S$-algebraic space finite flat of finite presentation over $S$ of constant degree $\delta_i$, $i=1,\cdots, n$. Suppose that $X$ has sections over all $S_i$, $i=1,\cdots, n$.}
\smallskip

\emph{Then $X$ is a trivial $A$-torsor.}

\begin{proof} By considering $A$ as $\mathrm{Pic}^o_{A^{*}/S}$, $A^{*}$ being the dual abelian algebraic space of $A$, one defines for each $i=1,\cdots, n$ the norm homomorphism
\[N_i: \prod_{S_i/S} A_{S_i}\to A.\] The composition
\[A\to \prod_{S_i/S}A_{S_i}\stackrel{N_i}{\longrightarrow} A\] is equal to $\delta_i\mathrm{Id}_A$, where
\[A\to \prod_{S_i/S}A_{S_i},\ x\mapsto x_{S_i}=x\times_SS_i\] is the adjunction morphism associated with the pair of adjoint functors
\[\prod_{S_i/S}-\ ,\ S_i\times_S-.\]

Let $\sigma_i$ be a section of $X$ over $S_i$, $i=1,\cdots, n$, and choose integers $e_1,\cdots, e_n$ such that $e_1\delta_1+\cdots+e_n\delta_n=1$. Consider the $S$-morphism
\[q: X\to A,\ x\mapsto \sum^n_{i=1}e_i.N_i(x_{S_i}-\sigma_i),\] where $x_{S_i}-\sigma_i$ denotes the unique local section $a_i$ of $A_{S_i}$ satisfying $a_i+\sigma_i=x_{S_i}$. For each local $S$-section $a$ (resp. $x$) of $A$ (resp. $X$), one has
\[q(a+x)=\sum^n_{i=1}e_i(N_i(a_{S_i})+N_i(x_{S_i}-\sigma_i))=(\sum^n_{i=1}e_i\delta_i.a)+q(x)=a+q(x).\]

Being thus an $A$-equivariant morphism between $A$-torsors, $q$ is an isomorphism.

\end{proof}

\smallskip

{\bf Theorem 4.9.} --- \emph{Keep the assumptions of $(4.1)$. Assume furthermore that $X_{\overline{s}}$ is connected, proper, of total multiplicity prime to the characteristic of $k(\overline{s})$ and that $X$ is regular.}
\smallskip

\emph{Then there is a non-degenerate abelian fibration $F/E$ over an $S$-algebraic stack $E$ with $E\times_St=t$ which extends $X_t$ over $t$ where $E$ is finite flat over $S$, tame along $s$ and regular. The identity $X_t=F\times_EE_t$ extends to a proper $S$-morphism $p$ from $X$ to $F$ which is \'{e}tale at precisely the points outside the image of the uniruled irreducible components of $X_{\overline{s}}$. Such $(F/E, p)$ is unique up to unique $S$-isomorphisms and its formation commutes with every formally smooth faithfully flat base change $T\to S$ of spectra of discrete valuation rings.}

\smallskip

{\bf Theorem 4.10.} --- \emph{Keep the hypothesis of $(4.9)$ except that $X$ be regular. Assume that $X$ is at each of its geometric points either geometrically factorial of equal characteristic or regular. Assume furthermore that $X_{\overline{s}}$ does not have uniruled irreducible components.}
\smallskip

\emph{Then $X$ is regular and is a universal almost non-degenerate $S$-abelian fibration and $X_{\overline{s}}$ does not contain $\overline{s}$-rational curves. If $A/E$ denotes the albanese of $X/S$, the action of $A_t$ on $X_t$ extends uniquely to an action on $X$ by the $S$-N\'{e}ron model
\[\prod_{E/S}A\] of $A_t$.}

\smallskip

4.11. \emph{Proof of} (4.9).
\smallskip

--- \emph{Case where $S$ is strictly henselian }:
\smallskip

Let $\overline{t}$ be the spectrum of a separable closure of $k(t)$ and let $\overline{\eta}$ be a geometric generic point of $X_{\overline{t}}$. The projection $X_{\overline{t}}\to X$ induces the specialization homomorphism (SGA 1 X 2)
\[sp: \pi_1(X_{\overline{t}}, \overline{\eta})\to \pi_1(X, \overline{\eta}).\] The image of $sp$ is a normal subgroup of finite index (\cite{raynaud_specialization} 6.3.5) and its associated monodromy representation
\[\pi_1(X, \overline{\eta})\to \mathrm{Coker}(sp)=G\] corresponds by Galois theory to an $X$-algebraic space $X'$ which is connected and finite \'{e}tale Galois over $X$ with Galois group $G$. 
\smallskip

Let $S'=\mathrm{Spec}\ \Gamma(X', \mathcal{O}_{X'})$, which is the spectra of a discrete valuation ring. Let the generic (resp. closed) point of $S'$ be $t'$ (resp. $s'$). Then $X'_{t'}=X_t\times_tt'$. By \cite{raynaud_specialization} 6.3.5+6.3.7, $X'_{s'}$ is of total multiplicity $1$, $k(s)=k(s')$ and $G$ is cyclic of order equal to the total multiplicity $\delta$ of $X_s$, as $\delta$ is prime to the characteristic of $k(s)$ and $X$ being regular has geometrically factorial local rings. 
\smallskip

Let the maximal points of $X'_{s'}$ be $x'_1,\cdots, x'_n$ and let $Z'_i$ be the closed image of $\mathrm{Spec}(\mathcal{O}_{X'_{s'}, x'_i})\to X'_{s'}$ in $X'_{s'}$, $i=1,\cdots, n$. By \cite{raynaud_specialization} 7.1.2, for each $i=1,\cdots, n$, there is a regular closed $S'$-immersion $S'_i\hookrightarrow X'$ such that $S'_i$ is finite flat over $S'$ of rank equal to the total multiplicity $\delta'_i$ of $Z'_i/s'$ and such that $S'_i$ intersects $X'_{s'}$ at one unique point of $Z'_i\backslash\sum_{j\neq i}Z'_j$. The greatest common divisor of $\delta'_1,\cdots, \delta'_n$ is by definition the total multiplicity of $X'_{s'}$, that is, $1$. Thus by (4.8) the $(A_t\times_tt')$-torsor $X'_{t'}=X_t\times_tt'$ admits a $t'$-point; by means of one such $t'$-point, one identifies $X'_{t'}$ with $A_t\times_tt'$.
\smallskip

Write $W'$ for the open sub-algebraic space of $X'$ which consists of all points at which $X'\to S'$ is smooth. Each $t'$-point of $X'_{t'}$ uniquely extends to an $S'$-section of $W'$, as $X'$ is proper over $S'$ and regular. That is, $W'$ is a weak $S'$-N\'{e}ron model of $W'_{t'}=X'_{t'}$ (\cite{neron_model} 3.5/1). Let $d=\mathrm{dim}(W'_{t'})$. If one chooses a non-zero section $\omega'\in\Gamma(W', \Omega^d_{W'/S'})$ such that the divisor $\mathrm{Div}_{W'}(\omega')$ has support strictly contained in $W'_{s'}$, the open $U'=W'-\mathrm{Supp}(\mathrm{Div}_{W'}(\omega'))$ has an $S'$-birational group law which extends the group structure of $U'_{t'}=X'_{t'}=A_t\times_tt'$ over $t'$. 
\smallskip

One argues as in (4.7) that the $S'$-N\'{e}ron model of $A_t\times_tt'$ is an $S'$-abelian scheme $A'$ and that the regular $S'$-minimal model $F'$ of $X'_{t'}$ is a trivial $A'$-torsor. 
The identity $X'_{t'}=F'_{t'}$ extends uniquely by (2.1) to an $S'$-morphism $p': X'\to F'$, which is equivariant with respect to the canonical action of $G$ on $X'$ and on $F'$. Moreover, $p'$ is proper surjective and is \'{e}tale at precisely the points of $U'$ and $U'$ is the complement of the image of the uniruled irreducible components of $X'_{\overline{s'}}$, where $\overline{s'}$ denotes the spectrum of an algebraic closure of $k(s')$.
The quotient of $p'$ by $G$, 
\[p=[p'/G]: [X'/G]=X\to [F'/G]=F,\] is consequently proper and is \'{e}tale at precisely the points outside the image of the uniruled irreducible components of $X_{\overline{s}}$. The projection
\[[F'/G]=F\to [S'/G]=E\] is a non-degenerate abelian fibration with albanese $[A'/G]=A$. The algebraic stack $E$ is finite flat over $S$, regular and tame along $s$, as $S'$ is. Over $t$, one has
$E\times_St=[S'_{t}/G]=[t'/G]=t$.
\smallskip

The formation of $(F/E, p)$ evidently commutes with every formally smooth faithfully flat base change $T\to S$ of spectra of strictly henselian discrete valuation rings.
\smallskip

It remains to characterize $(F/E, p)$ up to unique $S$-isomorphisms. Let $(F^{\natural}/E^{\natural}, p^{\natural})$ be an alternative with albanese $A^{\natural}$. Let $U=U'/G$ be the complement in $X$ of the image of the uniruled irreducible components of $X_{\overline{s}}$. There is a unique $U$-isomorphism of $U$-abelian algebraic spaces $A^{\natural}\times_{E^{\natural}}U=A\times_EU$ extending the identity morphism of $A_t\times_tX_t$. For, the restriction functor from the category of $U$-abelian algebraic spaces to the category of $X_t$-abelian algebraic spaces is fully faithful. Let $U_1^{\natural}$ be the open sub-$U$-algebraic space of $A^{\natural}\times_{E^{\natural}}U=A\times_EU=A_U$ which has image $U\times_{E^{\natural}}U$ by the isomorphism
\[r^{\natural}=(\mu^{\natural}, p_2^{\natural}): A^{\natural}\times_{E^{\natural}}F^{\natural}\ \widetilde{\to}\ F^{\natural}\times_{E^{\natural}}F^{\natural},\] where $p_2^{\natural}$ (resp. $\mu^{\natural}$) is the second projection (resp. represents the action of $A^{\natural}$ on $F^{\natural}$). Write $U_1^{\natural'}$ for the inverse image of $U_1^{\natural}$ by the projection $A'_{U'}\to A_U$ and write $x'$ (resp. $y'$) for the generic point of $U'_{s'}$ (resp. $A'_{U'}\times_{S'}s'$). Now $U^{\natural'}$ contains $y'$, and $r^{\natural}$ induces an $S$-morphism 
\[(\mathrm{Spec}(\mathcal{O}_{U^{\natural'}_1, y'})\rightrightarrows\mathrm{Spec}(\mathcal{O}_{U', x'}))\to (U\times_{E^{\natural}}U\rightrightarrows U),\] which in turn by quotient induces an $S$-morphism
$S'\to E^{\natural}$. This latter is smooth surjective, since the composition
$U'\to S'\to E^{\natural}$ is. As both $S'$ and $E^{\natural}$ are finite over $S$, $S'\to E^{\natural}$ is finite \'{e}tale surjective and hence, for each integer $n\geq 0$, $\mathrm{cosq}_o(S'/E^{\natural})_n$ is the normalization of $S$ in $\mathrm{cosq}_o(t'/t)_n$. So $E^{\natural}=E$. To prove that $A^{\natural}=A$, it suffices to prove the equality
$A^{\natural}\times_ES'=A\times_ES'=A'$ and that the descent data on $A'$ relative to $S'\to E$ corresponding to $A^{\natural}$ and to $A$ coincide. One has a unique $S'$-isomorphism of $S'$-abelian algebraic spaces $A^{\natural}\times_ES'=A\times_ES'$ which extends the identity morphism of $A'_{t'}$, for the restriction functor from the category of $S'$-abelian algebraic spaces to the category of $t'$-abelian algebraic spaces is fully faithful. The coincidence of the two descent data is deduced in the same way, since the restriction functor from the category of abelian algebraic spaces over $G_{S'}$ (resp. $(G\times G)_{S'}$) to the category of abelian algebraic spaces over $G_{t'}$ (resp. $(G\times G)_{t'}$) is fully faithful. The identity $F^{\natural}=F$ follows by a similar argument based on the uniqueness of regular minimal models. Finally, $p^{\natural}=p$, as one has $(p^{\natural}\times_ES')|X'_{t'}=(p\times_ES')|X'_{t'}$.
\smallskip

--- \emph{General case }:
\smallskip

Let $S_{(\overline{s})}$ denote the strictly henselization of $S$ at $\overline{s}$. Let $\pi\in\Gamma(S, \mathcal{O}_S)$ be a uniformizer and $f: X\to S$ the structural morphism. Notice that the cycle $\Delta=f^*\mathrm{Div}_S(\pi)/\delta$ is integral, where $\delta$ denotes the total multiplicity of $X_{\overline{s}}$, for $\Delta\times_SS_{(\overline{s})}$ is on $X\times_SS_{(\overline{s})}$. With $\Delta$ one associates a canonical $\mu_{\delta}$-torsor on $X$ for the \'{e}tale topology, $X'\to X$, which after the base change $S_{(\overline{s})}\to S$ corresponds to the specialization homomorphism of fundamental groups above. In particular, it suffices to define $E$ to be $[S'/\mu_{\delta}]$, where $S'$ is defined to be $\mathrm{Spec}\ \Gamma(X', \mathcal{O}_{X'})$. There is by quotient by $\mu_{\delta}$ an $S$-morphism
\[[X'/\mu_{\delta}]=X\to [S'/\mu_{\delta}]=E.\] Let $t'$ be the generic point of $S'$. One verifies after the base change $S_{(\overline{s})}\to S$ that the $S'$-N\'{e}ron model of $A_t\times_tt'$ is an $S'$-abelian scheme $A'$ and that the regular $S'$-minimal model of $X'_{t'}=X_t\times_tt'$ is an $A'$-torsor $F'$ on $S'$ for the \'{e}tale topology. On $F'\to S'$, $\mu_{\delta}$ acts compatibly. The non-degenerate abelian fibration
\[[F'/\mu_{\delta}]=F\to [S'/\mu_{\delta}]=E\] has albanese $[A'/\mu_{\delta}]=A$. The identity $X'_{t'}=F'_{t'}$ extends by (2.1) to a unique $S'$-morphism $p': X'\to F'$. Write $p=[p'/\mu_{\delta}]$, which is proper and is \'{e}tale at precisely the points outside the image of the uniruled irreducible components of $X_{\overline{s}}$. This $(F/E, p)$ is unique up to unique $S$-isomorphisms and its formation commutes with every formally smooth faithfully flat base change $T\to S$ of spectra of discrete valuation rings, as one verifies after the base change $S_{(\overline{s})}\to S$.

\begin{flushright}
$\square$
\end{flushright}

4.12. \emph{Proof of} (4.10).
\smallskip

Keep the notations of (4.11). As $X_{\overline{s}}$ by hypothesis does not have uniruled irreducible components, the morphism $p': X'\to F'$ is thus an isomorphism (4.2) and $X'_{\overline{s'}}$, hence $X_{\overline{s}}$ as well, does not contain rational curves. Clearly, $X=F$ is an almost non-degenerate $S$-abelian fibration with ramification stack $E$ and albanese $A$. This fibration is universal. For, as $\delta$ is prime to the residue characteristics of $S$, the formation of the quotient $S'/\mu_{\delta}$ commutes with every base change $T\to S$. 
\smallskip

Write
\[\overline{A}=\prod_{E/S}A,\] which is the kernel of the diagram
\[(d_o^*, d_1^*): \prod_{S'/S}A'\rightrightarrows \prod_{S''/S}A'',\] where $S''=\mu_{\delta}\times_SS'$, $d_1, d_o: S''\to S'$ respectively denotes the second projection and represents the action of $\mu_{\delta}$ on $S'$, and $A''=A\times_ES''=d_o^*A'=d_1^*A'$. In particular, $\overline{A}$ is a separated $S$-group scheme of finite type, since both
\[\prod_{S'/S}A',\ \prod_{S''/S}A''\] are separated $S$-group schemes of finite type (\cite{neron_model} 7.6/4). Moreover, $\overline{A}(t^{hs})=\overline{A}(S_{(\overline{s})})$, where $t^{hs}$ denotes the generic point of $S_{(\overline{s})}$, as $A'$ is an $S'$-abelian scheme. Next, let $T$ be an affine $S$-scheme and let $T_o$ be a closed sub-$S$-scheme of $T$ defined by an ideal $I$ with $I^2=0$. By applying the functor $H^0(\mu_{\delta}, -)$ to the exact sequence of $\mu_{\delta}$-modules
\[0\to\Gamma(T\times_SS', I\otimes_S\mathrm{Lie}(A'/S'))\to A'(T\times_SS')\to A'(T_o\times_SS')\to 0\] one obtains a surjection ($\delta$ invertible on $S$)
\[\overline{A}(T)=H^0(\mu_{\delta}, A'(T\times_SS'))\to \overline{A}(T_o)=H^0(\mu_{\delta}, A'(T_o\times_SS')).\] This shows that $\overline{A}$ is formally smooth over $S$. So $\overline{A}$ is the $S$-N\'{e}ron model of $\overline{A}_t=A_t$. 
\smallskip

Write the action of $A_t$ on $X_t$ as $\mu_t: A_t\times_tX_t\to X_t$. By (2.1) $\mu_t$ uniquely extends to an $S$-morphism $\mu: \overline{A}\times_SX\to X$, as $\overline{A}\times_SX$ is regular connected and as $X_{\overline{s}}$ does not contain $\overline{s}$-rational curves. The $S$-binary law $\mu$ is associative and hence represents an action of $\overline{A}$ on $X$, as $A_t\times_tA_t\times_tX_t$ is dense in $\overline{A}\times_S\overline{A}\times_SX$ and $X$ is $S$-separated.

\begin{flushright}
$\square$
\end{flushright}

\smallskip

5. \emph{Non-uniruled abelian fibrations in characteristic zero. Purity.}
\smallskip

{\bf Lemma 5.1.} --- \emph{Let $S$ be a locally noetherian algebraic space and $U$ an open sub-algebraic space of $S$ with $\mathrm{prof}_{S-U}(S)\geq 2$.}
\smallskip

\emph{Then the functor $A\mapsto A|U$, from the category of $S$-abelian algebraic spaces to the category of $U$-abelian algebraic spaces, is fully faithful. It is an equivalence if $S$ is normal of residue characteristics zero and pure along $S-U$} (SGA 2 X 3.1)\emph{, in particular, if $S$ is regular of residue characteristics zero.}

\begin{proof} The full-faithfulness of the functor $A\mapsto A|U$ follows by \cite{lemme de Gabber} Proposition (1), 3). The assertion on the equivalence is \cite{grothendieck_abelian} 4.2+4.5. 

\end{proof}

\smallskip

{\bf Lemma 5.2.} --- \emph{Let $S$ be a noetherian local scheme with closed point $s$. Let $U=S-\{s\}$. Let $A$ be an $S$-abelian algebraic space with structural morphism $f$ and zero section $e$. Let $A^*=\mathrm{Pic}^o_{A/S}$ be the dual $S$-abelian algebraic space of $A$ with structural morphism $f^*$ and zero section $e^*$.}
\smallskip

\emph{Then the two following statements hold when $A$ is para-factorial} (SGA 2 XI 3.1)\emph{ along $f^{-1}(s)$ }:
\smallskip

1) \emph{Each $U$-section of $f^*|U$ extends uniquely to an $S$-section of $f^*$.}
\smallskip

2) \emph{Each $f|U$-fiberwise numerically trivial invertible module on $f^{-1}(U)$ rigidified along $e|U$ extends uniquely to an $f$-fiberwise numerically trivial invertible module on $A$ rigidified along $e$.}

\begin{proof} Note that these two are the same statements. By the hypothesis that $A$ is para-factorial along $f^{-1}(s)$, one has $\mathrm{prof}_s(S)\geq 2$ and that each invertible module $L$ on $f^{-1}(U)$ extends up to unique isomorphisms to a unique invertible module $\overline{L}$ on $A$. It is evident that $\overline{L}$ is $f$-fiberwise numerically trivial (resp. has a unique rigidification along $e$ extending that of $L$ along $e|U$) if $L$ is $f|U$-fiberwise numerically trivial (resp. rigidified along $e|U$).  

\end{proof}

\smallskip

{\bf Lemma 5.3.} --- \emph{Let $S$ be a noetherian local scheme with closed point $s$. Let $X$ be an $S$-smooth algebraic space with structural morphism $f$. Then $X$ is para-factorial along $f^{-1}(s)$ in the following cases }:
\smallskip

i) \emph{$S$ is regular of dimension $\geq 2$.}
\smallskip

ii) \emph{$S$ is of equal characteristic and geometrically para-factorial at $s$.}

\begin{proof} Let $\overline{s}$ be the spectrum of a separable closure of $k(s)$ and $S_{(\overline{s})}$ the strict localization of $S$ at $\overline{s}$. Recall that one says that $S$ is \emph{geometrically para-factorial} at $s$ if $S_{(\overline{s})}$ is para-factorial along $\overline{s}$. Case i) is classical and ii) is \cite{boutot} III 2.14.

\end{proof}
\smallskip

{\bf Definition 5.4.} --- \emph{Let $S$ be an algebraic space and $U$ an open sub-algebraic space of $S$. We say that $S$ is $A$-pure along $S-U$ if for every smooth morphism $S'\to S$ the functor $A\mapsto A|U'$, from the category of $S'$-abelian algebraic spaces to the category of $U'$-abelian algebraic spaces, is an equivalence, where $U'=U\times_SS'$. We say that $S$ is strictly $A$-pure along $S-U$ if furthermore for every smooth morphism $S'\to S$ with $U'=U\times_SS'$ and every $S'$-abelian algebraic space $A$, each $U'$-section of $A$ extends uniquely to an $S'$-section of $A$. We say that $S$ is $A$-pure (resp. strictly $A$-pure) at a geometric point $s$ if its strict localization $S_{(s)}$ at $s$ is $A$-pure (resp. strictly $A$-pure) along $s'$, where $s'$ is the closed point of $S_{(s)}$.}

\smallskip

{\bf Example 5.5.} --- By (5.1) and by SGA 4 XV 2.1, if an algebraic space $S$ is of residue characteristics zero locally noetherian normal and \emph{pure} (SGA 2 X 3.2) at a geometric point $s$, then $S$ is $A$-pure at $s$. If furthermore the strict localization $S_{(s)}$ is para-factorial along its closed point, then by (5.2)+(5.3) $S$ is strictly $A$-pure at $s$. 
\smallskip

{\bf Lemma 5.6.} --- \emph{Let $S$ be an algebraic space and $U$ an open sub-algebraic space of $S$ such that $S$ is strictly $A$-pure along $S-U$. For $i=1, 2$, let $A_i$ be an $S$-abelian algebraic space and $X_i$ an $A_i$-torsor on $S$ for the \'{e}tale topology.}
\smallskip

\emph{Then each $U$-morphism from $X_1|U$ to $X_2|U$ extends uniquely to an $S$-morphism from $X_1$ to $X_2$.}

\begin{proof} The question being an \'{e}tale local question on $S$, one may assume the torsors $X_i$, $i=1, 2$, to be trivial. 
\smallskip

Each $U$-morphism $q: A_1|U\to A_2|U$ is the unique composite of a translation ($a_2\mapsto a_2+q(0)$) and a $U$-group homomorphism $p: A_1|U\to A_2|U$ (``Geometric Invariant Theory'' 6.4). As by hypothesis $S$ is strictly $A$-pure along $S-U$, the $U$-section $q(0)=\sigma$ and the $U$-group homomorphism $p$ extend uniquely to an $S$-section $\overline{\sigma}$ and an $S$-group homomorphism $\overline{p}$, hence the claim. 

\end{proof}

\smallskip

{\bf Proposition 5.7.} --- \emph{Let $S$ be an algebraic space and $U\to X$ an open immersion of $S$-algebraic spaces such that $X$ is strictly $A$-pure along $X-U$. Assume that on $U$ there is given an almost non-degenerate $S$-abelian fibration structure with defining $S$-groupoid $U_.$.}
\smallskip

\emph{Then up to unique isomorphisms there exists a unique almost non-degenerate $S$-abelian fibration structure on $X$ with $S$-groupoid $X_.$ such that $d_1: X_1\to X_o=X$ restricts to $d_1: U_1\to U_o=U$ on $U$.}

\begin{proof} As $X$ is $A$-pure along $X-U$, there exist unique cartesian diagrams in the category of $S$-algebraic spaces :
{\[\xymatrix{
                  U_1 \ar[r]^{j'} \ar[d]_{d_1}  & A' \ar[d]^{f'} \\
                  U \ar[r]_{}                              & X}\] }
{\[\xymatrix{
                  U_1 \ar[r]^{j} \ar[d]_{d_o}   &  A\ar[d]^{f} \\
                  U \ar[r]_{}                              & X}\] }whose vertical arrows have abelian algebraic space structures. Let $p'$ (resp. $p$) denote the projection of $A'\times_XA$ onto $A'$ (resp. $A$). The diagonal immersion
\[(j', j): U_1\to A'\times_XA\] satisfies 
\[p'(j', j)=j',\ p(j', j)=j.\] 

--- \emph{There is a unique section $i'$ of the abelian algebraic space structure $p'$ such that $i'j'=(j', j)$.}
\smallskip

Indeed, as $j'$ is the base change of $U\hookrightarrow X$ by the smooth morphism $f'$, $A'$ is strictly $A$-pure along $A'-j'(U_1)$. So $(j', j)$, considered as a $U_1$-section of $p'$, extends uniquely to a section $i'$ of $p'$. That is, $p'i'=\mathrm{Id}_{A'}$ and $i'j'=(j', j)$.
\smallskip

--- \emph{The following diagrams are commutative and cartesian }:
{\[\xymatrix{
                  U_1 \ar[r]^{j'} \ar[d]_{d_1}  & A' \ar[d]^{f'} \\
                  U \ar[r]_{}                              & X}\] }
{\[\xymatrix{
                  U_1 \ar[r]^{j'} \ar[d]_{d_o}  & A' \ar[d]^{fpi'} \\
                  U \ar[r]_{}                              & X}\] }which, from now on, we rewrite as :                   
{\[\xymatrix{
                  U_1 \ar[r]^{} \ar[d]_{d_1}    & X_1 \ar[d]^{d_1} \\
                  U_o \ar[r]_{}                          & X_o}\] }                  
{\[\xymatrix{
                  U_1 \ar[r]^{} \ar[d]_{d_o}    & X_1 \ar[d]^{d_o} \\
                  U_o \ar[r]_{}                          & X_o}\] }Next, form the cartesian diagram :
{\[\xymatrix{
                  X_2 \ar[r]^{d_o} \ar[d]_{d_2}  & X_1 \ar[d]^{d_1} \\
                  X_1 \ar[r]_{d_o}                       & X_o}\] }

--- \emph{There is a unique morphism $d_1: X_2\to X_1$ such that the following diagram commutes and is cartesian }:
{\[\xymatrix{
                  X_2 \ar[r]^{d_1} \ar[d]_{d_2}  & X_1 \ar[d]^{d_1} \\
                  X_1 \ar[r]_{d_1}                       & X_o}\] }
                  
Indeed, as $U_1\hookrightarrow X_1$ is the base change of $U\hookrightarrow X$ by the smooth morphism $d_1$, $X_1$ is strictly $A$-pure along $X_1-U_1$. By (5.6) the cartesian diagram (SGA 3 V 1)
{\[\xymatrix{
                  U_2 \ar[r]^{d_1} \ar[d]_{d_2}  & U_1 \ar[d]^{d_1} \\
                  U_1 \ar[r]_{d_1}                       & U_o}\] }whose vertical arrows have abelian algebraic space structures has a unique extension as above claimed, for one verifies that :
\smallskip

\noindent i) \emph{The base change of $d_2: X_2\to X_1$ by $U_1\hookrightarrow X_1$ is $d_2: U_2\to U_1$.}
\smallskip

\noindent ii) \emph{The base change of $d_1: X_1\to X_o$ by $U_1\hookrightarrow X_1\stackrel{d_1}{\longrightarrow}X_o$ is the base change of $d_1: U_1\to U_o$ by $d_1: U_1\to U_o$.}
\smallskip

--- \emph{The above $d_1: X_2\to X_1$ fits into the following diagram which is commutative and cartesian }:                  
{\[\xymatrix{
                  X_2 \ar[r]^{d_1} \ar[d]_{d_o}  & X_1 \ar[d]^{d_o} \\
                  X_1 \ar[r]_{d_o}                       & X_o}\] }For, $X_1$ is strictly $A$-pure along $X_1-U_1$ and one has similarly the cartesian diagram :
{\[\xymatrix{
                  U_2 \ar[r]^{d_1} \ar[d]_{d_o}  & U_1 \ar[d]^{d_o} \\
                  U_1 \ar[r]_{d_o}                       & U_o}\] } It is now immediate that one has obtained the desired $S$-groupoid $X_.$ (cf. SGA 3 V 1).
                  
\end{proof}

\smallskip

{\bf Proposition 5.8.} --- \emph{Let $S$ be a noetherian normal integral scheme, $t$ the generic point of $S$, $A_t$ a $t$-abelian variety and $X_t$ an $A_t$-torsor on $t$ for the \'{e}tale topology. Assume that, for each strict henselization $S'$ of $S$ at a geometric codimension $1$ point $s$, if $t'$ (resp. $s'$) denotes the generic (resp. closed) point of $S'$, $X_t\times_tt'$ extends to a separated $S'$-algebraic space $X'$ of finite type such that $X'$ is normal integral and at each of its geometric codimension $\geq 2$ points either regular or pure geometrically para-factorial of equal characteristic and that $X'_{s'}$ is non-empty, separable, proper and does not have uniruled irreducible components.}
\smallskip

\emph{Then, if $S$ is $A$-pure at all its geometric points of codimension $\geq 2$, there exists up to unique isomorphisms a unique $S$-abelian algebraic space $A$ extending $A_t$.}

\begin{proof} Recall that the formation of N\'{e}ron models commutes with strict localization. Thus, the N\'{e}ron model of $A_t$ at every codimension $1$ point of $S$ is by (4.3) an abelian scheme. So, as $S$ is $A$-pure at all its geometric points of codimension $\geq 2$, there is up to unique isomorphisms a unique extension of $A_t$ to an $S$-abelian algebraic space $A$.

\end{proof}

\smallskip

{\bf Proposition 5.9.} --- \emph{Keep the notations of $(5.8)$. Assume that $S$ is of residue characteristics zero pure at all its points of codimension $\geq 2$ and that there is an open sub-scheme $R$ of $S$ which consists precisely of all points of $S$ where $S$ is regular.}
\smallskip

\emph{Then $X_t$ extends to an $A$-torsor $X$ on $S$ for the \'{e}tale topology. Such an extension is unique up to unique isomorphisms if $S$ is geometrically para-factorial along $S-R$.}

\begin{proof} Note that $S$ is by (5.5) $A$-pure at all its geometric points of codimension $\geq 2$. So (5.8) applies. 
\smallskip

As the formation of regular minimal models commutes with strict localization, the regular minimal model of $X_t$ at each codimension $1$ point $s$ of $S$ is by (4.3) a torsor for the \'{e}tale topology under the localization of $A$ at $s$. As $R$ is strictly $A$-pure at all its points of codimension $\geq 2$, there exist, by (5.6) and a ``passage \`{a} la limite'', an open sub-scheme $V$ of $R$ with $\mathrm{codim}(R-V, R)\geq 2$ and an $A|V$-torsor $Z$ on $V$ for the \'{e}tale topology such that $Z$ extends $X_t$. 
\smallskip

This torsor $Z$ is by \cite{raynaud_thesis} XIII 2.8 iv) of finite order. Namely, there exist an integer $n\geq 1$ and an ${}_nA|V$-torsor $P$ on $V$ for the \'{e}tale topology such that 
\[Z=P\stackrel{{}_nA|V}{\wedge}A|V,\] where ${}_nA=\mathrm{Ker}(n.\mathrm{Id}_A)$. 
\smallskip

As $S$ is pure at all its points of codimension $\geq 2$, thus in particular pure along $S-V$, there is a unique finite \'{e}tale $S$-scheme $\overline{P}$ which restricts to $P$ on $V$. By the purity of $S$ along $S-V$ again, $\overline{P}$ is in a unique way an ${}_nA$-torsor on $S$ for the \'{e}tale topology and hence
\[X=\overline{P}\stackrel{{}_nA}{\wedge}A\] extends $X_t$. Such an extension is by (5.6) unique up to unique isomorphisms if $S$ is geometrically para-factorial at all its points of codimension $\geq 2$, or equivalently, at all points of $S-R$.

\end{proof}

\smallskip

{\bf Theorem 5.10.} --- \emph{Let $S$ be an integral scheme with generic point $t$ and $X$ an $S$-algebraic space with structural morphism $f$. Assume that $X$ is locally noetherian normal integral of residue characteristics zero and at all its geometric codimension $\geq 2$ points pure and geometrically para-factorial. Assume furthermore that $f^{-1}(t)$ is a non-degenerate $t$-abelian fibration and that, for each geometric codimension $1$ point $\overline{x}$ of $X$, $f\times_SS_{(\overline{s})}$ is separated of finite type and flat at $\overline{x}$ and the geometric fiber $f^{-1}(\overline{s})$ is proper and does not have uniruled irreducible components, where $S_{(\overline{s})}$ denotes the strict henselization of $S$ at $\overline{s}=f(\overline{x})$.}
\smallskip

\emph{Then there exists a unique almost non-degenerate abelian fibration structure on $f$ extending that of $f_t$.}

\begin{proof} One applies (4.10), (5.5) and (5.7).

\end{proof}

\smallskip

{\bf Proposition 5.11.} --- \emph{Keep the notations of $(5.10)$. Let $(X_., d_., s_.)$ denote the $S$-groupoid of $X/S$. Consider the following conditions }:
\smallskip

1) \emph{$f$ is proper, $S$ is excellent regular.}
\smallskip

2) \emph{$f$ is proper, $S$ is locally noetherian normal and at each of its points satisfies the condition $(W)$ $(\mathrm{EGA\ IV}\ 21.12.8)$.}
\smallskip

\emph{Then, if $1)$ (resp. $2)$) holds, $S$ is the cokernel of $(d_o, d_1)$ in the full sub-category of the category of $S$-algebraic spaces consisting of the $S$-algebraic spaces (resp. $S$-schemes) which are $S$-separated and locally of finite type over $S$.}

\begin{proof} Let $Z$ be an $S$-separated algebraic space locally of finite type over $S$ and $p: X\to Z$ an $S$-morphism satisfying $pd_o=pd_1$.
As the $t$-groupoid $X_{.t}$ is simply connected, $p_t: X_t\to Z_t$ factors through a unique $t$-point, say $\sigma_t$, of $Z_t$. It amounts to showing that when 1) holds (resp. when 2) holds and $Z$ is a scheme) such a $t$-point uniquely extends to an $S$-section of $Z$.
\smallskip

Replacing $Z$ by the closed image of $\sigma_t$ in $Z$, one may assume that $Z$ is integral and birational over $S$. In case $1)$, as $p$ is dominant, $X$ normal and $S$ excellent, one may by replacing $Z$ by its normalization assume that $Z$ is normal.
\smallskip

As $f$ is proper and $Z$ is $S$-separated, $p$ is proper and hence surjective. It suffices to show that in case $1)$ (resp. $2)$ where $Z$ is a scheme) $Z$ is \'{e}tale over $S$ (resp. $Z\to S$ is a local isomorphism at every point of $Z$). For, being proper birational, $Z\to S$ is then an isomorphism.
\smallskip

When 1) holds (resp. when 2) holds and $Z$ is a scheme), it suffices by the theorem of \emph{purity of branch locus} (2.4) (resp. the theorem of \emph{purity of branch locus} of van der Waerden, EGA IV 21.12.12) to show that $Z$ is $S$-\'{e}tale at each geometric codimension $1$ point $\overline{z}$ of $Z$ (resp. $Z\to S$ is a local isomorphism at each codimension $1$ point $z$ of $Z$).
\smallskip

Now each geometric maxmal point $\overline{x}$ of $p^{-1}(\overline{z})$ (resp. each maximal point $x$ of $p^{-1}(z)$) is of codimension $\leq 1$ in $X$, and the image of $\overline{x}$ (resp. $x$) in $S$, which is also the image of $\overline{z}$ (resp. $z$), is of codimension $\leq 1$ in $S$ by hypothesis. 
The projection $Z\to S$ being proper birational is an isomorphism when localized at every codimension $\leq 1$ point of $S$ and in particular is \'{e}tale at $\overline{z}$ (resp. a local isomorphism at $z$). 

\end{proof}
\smallskip

5.12. \emph{Question }: 
\smallskip

In (5.11), does $\mathrm{Coker}(d_o, d_1)=S$ hold in the category of $S$-algebraic spaces? 
\smallskip


\bibliographystyle{amsplain}

\begin{thebibliography}{1}

\bibitem{neron_model}
S.~Bosch, W.~L\"{u}tkebohmert, M.~Raynaud.
\newblock N\'{e}ron models.
\newblock \emph{Ergebnisse der Mathematik und ihrer Grenzgebiete}, 1990.

\bibitem{boutot}
J.~F.~Boutot.
\newblock Sch\'{e}ma de Picard Local.
\newblock \emph{Lecture Notes in Mathematics}, 632, 1978.


\bibitem{dejong}
A.~J.~De Jong.
\newblock Smoothness, semi-stability and alterations.
\newblock \emph{Publications Math\'{e}matiques de l'IH\'{E}S}, 83, 1996.


\bibitem{lemme de Gabber}
P.~Deligne.
\newblock Le lemme de Gabber.
\newblock \emph{Ast\'{e}risque}, 127, 1985.






\bibitem{grothendieck_rational}
A.~Grothendieck.
\newblock  Techniques de construction en g\'{e}om\'{e}trie analytique. X. Construction de l'espace de Teichm\"{u}ller.
\newblock \emph{S\'{e}minaire Cartan}, 13, no. 2, 1960--1961.

\bibitem{grothendieck_abelian}
A.~Grothendieck.
\newblock Un th\'{e}or\`{e}me sur les homomorphismes de sch\'{e}mas ab\'{e}liens.
\newblock \emph{Invent. math}. 2, 1966.


\bibitem{geometric invariant theory}
\newblock D.~Mumford.
\newblock Geometric Invariant Theory. \emph{Ergebnisse der Mathematik und ihrer Grenzgebiete}, 1965.



\bibitem{raynaud_minimal}
M.~Raynaud.
\newblock Passage au quotient par une relation d'equivalence plate.
\newblock \emph{Proceedings of a conference on local fields}, NUFFIC Summer School held at Driebergen, Edited by T. A. Springer, 1966.

\bibitem{raynaud_specialization}
M.~Raynaud.
\newblock Sp\'{e}cialisation du foncteur de Picard.
\newblock \emph{Publications Math\'{e}matiques de l'IH\'{E}S}, 38, 1970.


\bibitem{raynaud_thesis}
M.~Raynaud.
\newblock Faisceaux amples sur les sch\'{e}mas en groupes et les espaces homog\`{e}nes.
\newblock \emph{Lecture Notes in Mathematics}, 119, 1971.

\bibitem{basic}
Y.~Zong.
\newblock Basic finite \'{e}tale equivalence relations. arxiv.org/abs/1512.00097.


\end{thebibliography}

\end{document}